\documentclass[12pt,english]{article}
\usepackage[T1]{fontenc}
\usepackage[latin9]{inputenc}
\usepackage[letterpaper]{geometry}
\usepackage{mathrsfs}
\usepackage{amsmath}
\usepackage{amsthm}
\usepackage{amssymb}
\usepackage{graphicx}
\usepackage{setspace}

\makeatletter
\numberwithin{equation}{section}
\numberwithin{figure}{section}
\theoremstyle{definition}
\newtheorem*{example*}{\protect\examplename}
\theoremstyle{plain}
\newtheorem*{lem*}{\protect\lemmaname}
\theoremstyle{plain}
\newtheorem*{thm*}{\protect\theoremname}

\@ifundefined{date}{}{\date{}}
\usepackage{pstricks}
\usepackage{changepage}

\makeatother

\usepackage{babel}
\providecommand{\examplename}{Example}
\providecommand{\lemmaname}{Lemma}
\providecommand{\theoremname}{Theorem}

\begin{document}
\begin{doublespace}
\noindent \begin{center}
\textbf{\Large{}A CONSTRUCTIVE VERSION OF THE}\\
\textbf{\Large{} SYLVESTER-GALLAI THEOREM}{\large{} \vspace{0.2cm}
}{\large\par}
\par\end{center}
\end{doublespace}

\noindent \begin{center}
MARK MANDELKERN{\small{}}\\
\par\end{center}

\begin{center}
{\small{}New Mexico State University, Las Cruces, New Mexico, USA}\\
\emph{\small{}e-mail: }{\small{}mandelkern@zianet.com}\\
\emph{\small{}web:}{\small{} www.zianet.com/mandelkern}{\small\par}
\par\end{center}

\vspace{1.7cm}

\begin{adjustwidth*}{2cm}{2cm}

\textbf{\small{}Abstract.}{\small{} The Sylvester-Gallai Theorem,
stated as a problem by James Joseph Sylvester in 1893, asserts that
for any finite, noncollinear set of points on a plane, there exists
a line passing through exactly two points of the set. First, it is
shown that for the real plane $\mathbb{R}^{2}$ the theorem is constructively
invalid. Then, a well-known classical proof is examined from a constructive
standpoint, locating the nonconstructivities. Finally, a constructive
version of the theorem is established for the plane $\mathbb{R}^{2}$;
this reveals the hidden constructive content of the classical theorem.
The constructive methods used are those proposed by Errett Bishop.}\footnote{ 

\hspace{0.3cm}\emph{Mathematics Subject Classification (MSC2010):}
primary: 51M04, secondary: 03F65. } 

\end{adjustwidth*}
\begin{quotation}
\vspace{0.1cm}
\end{quotation}
\begin{center}
\textbf{1. Introduction }\\
\par\end{center}

The Sylvester-Gallai Theorem states that for any finite, noncollinear
set of points on a plane, there exists a line passing through exactly
two points of the set. 

The history of this problem is itself problematic. The notion, that
Gallai was the first to prove the theorem, appears to stem from the
submission of the problem to the \emph{American Mathematical Monthly}
in 1943 by P. Erd\H{o}s {[}Erd43{]}, while unaware of the 1893 statement
of the problem by J. J. Sylvester {[}Syl93{]}. When a solution by
R. Steinberg {[}Ste44{]} was published in the\emph{ Monthly} in 1944,\emph{
}an \emph{Editorial Note} stated that Erd\H{o}s had ``enclosed with
the problem an outline of Grünwald's {[}Gallai's{]} solution''. This
reference to a pre-1944 proof by Gallai, albeit an unpublished outline,
appears to be the basis for the designation \emph{Sylvester-Gallai
Theorem. }This attribution recurs in a 1982 statement by Erd\H{o}s
{[}Erd82, p.208{]}, ``In 1933 . . . I told this problem to Gallai
who very soon found an ingenious proof\,''; Gallai's proof, it seems,
was not published. Along with Steinberg's 1944 solution, the \emph{Monthly}
noted that solutions were also received from R. C. Buck and N. E.
Steenrod. The earliest known published proof of the theorem, appearing
in 1941, is due to E. Melchior {[}Mel41{]}. In the present paper we
accede to common usage, refraining from use of the designation \emph{Sylvester-Melchior
Theorem}. 

There have been many different versions and proofs for this theorem.
V. Pambuccian {[}Pam09{]}, conducts reverse analyses of three proofs
of the theorem, leading to three different and incompatible axiom
systems. J. von Plato {[}Pla05{]} shows that the theorem holds intuitionistically
for sets of up to six points in a purely incidence-geometric setting,
and for up to seven points in an ordered geometric setting. See also
{[}KelMos58, Wil68, Cha70, Lin88, BorMos90, Chv04{]}.

We determine the constructive content of this theorem for the real
plane $\mathbb{R}^{2}$. First, we find that the theorem is constructively
invalid. Then we examine L. M. Kelly's 1948 proof,\footnote{Kelly's proof may be found in {[}Cox48{]} or {[}Cox61, pp. 65-66{]}. }
locating the nonconstructivities. Finally, adapting Kelly's method,
adding an hypothesis, and using strictly constructive methods, we
obtain a constructive version of the theorem. \\

\noindent \begin{center}
\textbf{2. Constructive methods} \\
\par\end{center}

The modern constructivist program began with L. E. J. Brouwer (1881-1966)
{[}Bro08{]}; recent work, using the strictest methods, follows the
work of Errett Bishop (1928-1983). A large portion of analysis is
constructivized by Bishop in \textit{Foundations of Constructive Analysis}
{[}B67{]}; this treatise also serves as a guide for constructive work
in other fields. This variety of constructivism does not form a separate
branch of mathematics, nor is it a branch of logic; it is intended
as an enhanced approach for all of mathematics. 

For the distinctive characteristics of Bishop-type constructivism,
as opposed to intuitionism or recursive function theory, see {[}BR87{]}.
Avoiding the \emph{Law of Excluded Middle }(LEM)\emph{, }constructive
mathematics is a generalization of classical mathematics, just as
group theory, a generalization of abelian group theory, avoids the
commutative law.\footnote{Constructive methods are described fully in {[}B67, B73, BB85{]};
see also {[}R82, BriMin84, M85, R99{]}. }

The initial phase of this program involves the rebuilding of classical
theories, using only constructive methods. The entire body of classical
mathematics is viewed as a wellspring of theories waiting to be constructivized. 
\begin{quotation}
\noindent Every theorem proved with {[}nonconstructive{]} methods
presents a challenge: to find a constructive version, and to give
it a constructive proof.

\noindent ~~~~~~~~~~~~~~~~~~~~~~~~~~~~~~~~~~~~~~~~~~~~~~~~~~~~~~~~~~~~~~~~~~~~-
Errett Bishop {[}B67, p. x{]} 
\end{quotation}
{\tiny{}~~~~~}{\tiny\par}

To clarify the constructive methods used here, we give examples of
familiar properties of the real numbers that are constructively \textit{invalid},
and also properties that are constructively \textit{valid}.\footnote{For more details, and other constructive properties of the real number
system, see {[}B67, BB85, BV06{]}. }

The following classical properties of a real number $\alpha$ are
constructively \emph{invalid}: 

\qquad{}(i)\emph{ Either $\alpha<0$ or $\alpha=0$ or $\alpha>0$}. 

\qquad{}(ii) \emph{If $\neg(\alpha=0),$ then $\alpha\ne0$.} 

Bishop constructs the real numbers using Cauchy sequences. Constructively,
the relation $\alpha\ne0$ does not refer to negation, but is given
a strong affirmative definition; one must construct an integer $n$
such that $1/n<|\alpha|$. 

Among the resulting constructively \emph{valid} properties of the
reals are the following: 

\qquad{}(a) \emph{For any real number $\alpha$, if $\neg(\alpha\ne0)$,
then $\alpha=0$. }

\qquad{}(b) \emph{For any real number $\alpha$, if $\neg(\alpha>0)$,
then $\alpha\le0$; if $\neg(\alpha<0)$, then $\alpha\ge0$. }

\qquad{}(c) \emph{Let }$\alpha$ and\emph{ $\beta$ be any real numbers,
with $\alpha<\beta$. For any real number $x$, either $x>\alpha$
or $x<\beta$.} 

\qquad{}(d) \emph{Let }$\alpha$ and\emph{ $\beta$ be any real numbers,
with $\alpha\neq\beta$. For any real number $x$, either $x\ne\alpha$
or $x\ne\beta$.} 

Property (c), known as the \emph{Constructive Dichotomy Principle,}
serves as a constructive substitute for the classical \emph{Trichotomy
Property}, which is constructively invalid. Property (d), which follows
from  (c), is called \emph{cotransitivity;} it is the classical contrapositive
of the transitive relation for equality. 

Points $A$ and $B$ on the real plane $\mathbb{R}^{2}$ are \emph{distinct},
written \emph{$A\neq B$}, if the distance between them is positive.
Points inherit properties from their coördinates; thus we have cotransitivity
for points: 

\qquad{}(e) \emph{Let $A$ and $B$ be any points on the real plane
$\mathbb{R}^{2}$, with $A\neq B$. For any point $X$, either $X\neq A$
or $X\neq B$. }

The condition \emph{$P$ lies on $l$, }written\emph{ $P\in l$,}
means that the distance $\rho(P,l)$, from the point $P$ to the line
$l$, is $0$, while \emph{$P$ lies outside $l$}, written\emph{
$P\notin l$,} means that the distance is positive.

The maximum and minimum of two real numbers are easily defined, using
Cauchy sequences. For any real number $\alpha$, we define $\alpha^{+}=\max\{\alpha,0\}$
and $\alpha^{-}=\max\{-\alpha,0\}$. Attempting to construct the maximum
or minimum of an arbitrary set of real numbers may, however, lead
to interesting complications. We say that a set is \emph{finite} if
its elements may be listed, $\{a_{1},a_{2},...,a_{n}\}$, with $n\geq1$,
but not necessarily distinctly.\footnote{Several different definitions and variations of the concept \emph{finite}
may be found in the constructive literature. The definition here is
close to the usual meaning of the term. } A finite set of real numbers has a minimum, but an arbitrary nonvoid
subset of a finite set of real numbers need not have a constructive
minimum. Even when a minimum does exist, it need not be attained by
an element in the set. These nonconstructivities will be established
in Section 5. \\

\noindent \begin{center}
\textbf{3. Constructively invalid statements} \\
\par\end{center}

Brouwerian counterexamples display the nonconstructivities in a classical
theory, indicating feasible directions for constructive work. To illustrate
the method, we give first an informal example on the plane $\mathbb{R}^{2}$.

\begin{example*}
\begin{flushleft}
If, for the real plane $\mathbb{R}^{2}$, there is a proof of the
statement
\par\end{flushleft}
\begin{flushleft}
\emph{\hspace{1cm}For any point $P$ and any line $l$, either $P$
lies on $l$, or $P$ lies outside $l$, }
\par\end{flushleft}
\begin{flushleft}
then we have a method that will either prove the \emph{Goldbach Conjecture,
}or construct a counterexample. 
\par\end{flushleft}
\end{example*}
\begin{proof}
Using a simple finite routine, construct a sequence $\{a_{n}\}_{n\geq2}$
such that $a_{n}=0$ if $2n$ is the sum of two primes, and $a_{n}=1$
if it is not. Now apply the statement in question to the point $P=(0,\Sigma a_{n}/n^{2})$,
with the $x$-axis as the line $l$. If $P$ lies on $l$, then we
have proved the Goldbach Conjecture, while if $P$ lies outside $l$,
then we have constructed a counterexample. 
\end{proof}
For this reason, such statements are said to be \emph{constructively
invalid.} If the Goldbach question is settled someday, then other
famous problems may still be ``solved'' in this way. The example
above will be applied in Section 5, locating one of the substantial
nonconstructivities in a classical proof of the Sylvester-Gallai Theorem. 

A \emph{Brouwerian counterexample} is a proof that a given statement
implies an omniscience principle. In turn, an \emph{omniscience principle}
would imply solutions or significant information for a large number
of well-known unsolved problems. This method was introduced by Brouwer
{[}Bro08{]} to demonstrate that reliance on the \emph{Law of Excluded
Middle }inhibits mathematics from attaining its full significance.
The omniscience principles are special cases of LEM. 

Omniscience principles are formulated in terms of binary sequences.
The zeros and ones may represent the results of a search for a solution
to a specific problem, as in the example above. These principles also
have equivalent statements in terms of real numbers. The following
omniscience principles will be used here: \\

\noindent \textbf{\textit{\emph{Limited principle of omniscience (LPO).
}}}\textit{For any binary sequence} $\{a_{n}\}$, \emph{either} $a_{n}=0$
\emph{for all} $n$, \emph{or there exists an integer} $n$ such that
$a_{n}=1$. Equivalently,\emph{ For any real number $\alpha$ with
$\alpha\geq0$, either $\alpha=0$ or $\alpha>0$. }\\

\noindent \textbf{\textit{\emph{Lesser limited principle of omniscience
(LLPO).}}}\emph{ }\textit{For any binary sequence} $\{a_{n}\}$, \emph{either
the first integer $n$} \emph{such that} $a_{n}=1$\emph{ (if one
exists) is even, or it is odd}. Equivalently,\emph{ For any real number
$\alpha$, either $\alpha\leq0$ or $\alpha\geq0$}.\emph{ }\\

A statement is considered \emph{constructively invalid} if it implies
an omniscience principle. The example above, with slight modification,
shows that the statement in question implies LPO; thus the statement
is constructively invalid. When omniscience principles are involved,
there is no need to mention a specific unsolved problem. LPO would
solve a great number of unsolved problems; LLPO would provide unlikely
partial solutions.\footnote{For more details concerning Brouwerian counterexamples, and other
omniscience principles, see ~{[}B73, R02, M89, BV06{]}. } \\

\noindent \begin{center}
\textbf{4. Brouwerian counterexample} \\
\par\end{center}

This will demonstrate that the Sylvester-Gallai Theorem is constructively
invalid for the plane $\mathbb{R}^{2}$.\\

\noindent \textbf{Example. }The statement,
\begin{quotation}
\noindent \emph{For any finite, noncollinear set of points $\mathscr{S}$
on the} \emph{real plane $\mathbb{R}^{2}$, there exists a line that
passes through exactly two points of $\mathscr{S}$,}
\end{quotation}
\noindent is constructively invalid; the statement implies LLPO. 

\includegraphics[viewport=-100bp 0bp 612bp 154bp,scale=0.5]{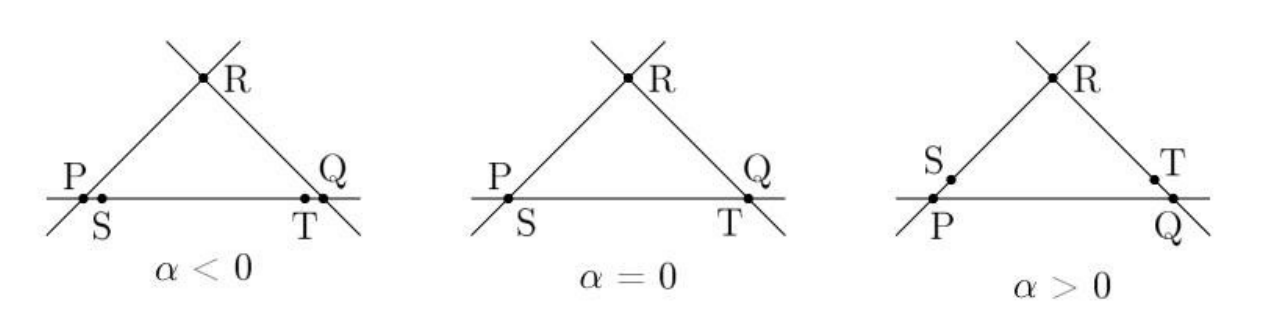}

\noindent \begin{adjustwidth*}{1.1cm}{1.6cm}

\noindent \textbf{Figure 1. }An omniscient view of the Brouwerian
counterexample. Constructively, we do not know whether the real number
$\alpha$ is $<0$, or $=0$, or $>0$; the classical \emph{Trichotomy
Property} is constructively invalid. Since we do not know about $\alpha$,
we cannot specify a line that passes through exactly two points. \\
\end{adjustwidth*}
\begin{proof}
Let $\alpha$ be any real number; either $|\alpha|>0$ or $|\alpha|<10^{-20}$.
Under the first condition, either $\alpha<0$ or $\alpha>0$, so we
already have the required conclusion; thus we may assume the second
condition. (Although the condition $|\alpha|<1/2$ would suffice to
keep the points in view, the condition chosen may bring to mind some
unsolved problem in number theory that has been checked a long way
out.) 

Define $\mathscr{S}=\{P,Q,R,S,T\}$, where $P=(-1,0)$, $Q=(1,0)$,
$R=(0,1)$, $S=(-1+|\alpha|,\alpha^{+})$, $T=(1-|\alpha|,\alpha^{+})$.
By hypothesis, there exists a line $l$ as specified in the statement;
there are at most ten possibilities for $l$. 

Consider first the case in which $l$ is one of the four lines\emph{
RP, RQ, RS, RT, }and suppose that $\alpha>0$. Then $\alpha^{+}=|\alpha|>0$,
so the point $S$ lies on the line $RP$, distinct from both points
$R$ and $P$. Similarly, the point $T$ lies on the line $RQ$. Thus
$l$ contains three distinct points of $\mathscr{S}$, a contradiction.
It follows that $\alpha\leq0$. 

Now consider the case in which $l$ is one of the six lines \emph{PQ,
PS, PT, QS, QT, ST,} and suppose that $\alpha<0$. Then $|\alpha|>0$,
while $\alpha^{+}=0$, so the points $S$ and $T$ lie on the line
$PQ$, distinct from each other and from $P$ and $Q$. Thus $l$
contains four distinct points of $\mathscr{S}$, a contradiction.
It follows that $\alpha\geq0$. 

Thus the statement implies LLPO. 
\end{proof}
\vspace{0.3cm}

\noindent \begin{center}
\textbf{5. Classical proof} \\
\par\end{center}

To obtain a constructive version of the Sylvester-Gallai Theorem,
we first subject Kelly's classical proof to examination from a constructive
viewpoint, locating the constructively invalid steps. Kelly's 1948
proof for the plane $\mathbb{R}^{2}$ is repeated below in brief form
(it appears in {[}Cox48{]} and {[}Cox61, pp. 65-66{]}). We follow
Bishop's suggestion, that theorems and proofs dependent on the \emph{Law
of Excluded Middle} be labeled as such. 
\begin{proof}
{[}LEM{]} (L. M. Kelly) Let $\mathscr{S}$ be a finite, noncollinear
set of points on the real plane $\mathbb{R}^{2}$, consider all the
lines that join two distinct points of $\mathscr{S}$, and consider
the distances from each point of $\mathscr{S}$ to each of these lines.
The set of all those distances that are positive has a minimum $d$,
attained by at least one line $l$ and one point $P$. Let $E$ denote
the foot of the perpendicular dropped from the point $P$ onto the
line $l$; the segment $PE$ then has length $d$. If there are three
distinct points of $\mathscr{S}$ that lie on $l$, then two of these
must lie on the same side of $E$. Denote these two points by $U$
and $V$, with $U$ closer to $E$ (possibly at $E$). Let $h$ be
the distance from the point $U$ to the line $PV$; then $h$ is one
of the positive distances considered in determining the minimum $d$.
However, it is apparent that $h$ is less than $d$, a contradiction.
Thus $l$ contains exactly two points of $\mathscr{S}$. 
\end{proof}
We encounter several constructive obstacles in this remarkable proof.
First, there is a problem in joining points to form the connecting
lines, since we cannot determine which pairs of given points are distinct;
considering distances, this would directly involve LPO. 

If we manage to clear the first obstacle, there is then the problem
of determining which points of $\mathscr{S}$ are at a positive distance
from which of the connecting lines; such a determination was shown
to be constructively invalid by the example in Section 3. Thus we
are unable to list the positive distances as a finite set of real
numbers, in order to construct a minimum. 

Although the set of positive distances may be a nonvoid subset of
a finite set of real numbers, it would be constructively invalid to
claim that such a set has a minimum. For a Brouwerian counterexample,
assume that such minimums exist, let $\alpha$ be a real number with
$\alpha\geq0$, and define $S=\{\alpha,1\}$. Either $\alpha>0$ or
$\alpha<1/2$; it suffices to assume the latter condition. Define
$T=\{x\in S\,:\,x>0\}$; then $T$ is a nonvoid subset of the finite
set $S$, so by hypothesis we may define $m=\min T$. Either $m>1/2$
or $m<1$. In the first case, $\neg(\alpha>0)$, so $\alpha=0$. In
the second case, since $m$ may be closely approximated by elements
of $T$, it follows that $\alpha\in T$, so $\alpha>0$. Thus LPO
results. In the proof above, if we cannot show that the set of positive
distances is finite, then we cannot construct a minimum. 

Assuming the resolution of the above difficulties, there remains a
problem with the selection of the special point-line pair. Even when
the minimum of a set of real numbers can be constructed, still we
cannot say that this minimum is attained by an element of the set.
For a Brouwerian counterexample, assume that minimums are attained,
let $\alpha$ be a real number, define $S=\{\alpha^{+},\alpha^{-}\}$,
and define $m=\min S$. By hypothesis, either $m=\alpha^{+}$, and
then $\neg(\alpha>0)$, so \emph{$\alpha\leq0$}, or $m=\alpha^{-}$,
and then $\neg(\alpha<0)$, so \emph{$\alpha\geq0$. }Thus LLPO results.\emph{
}Thus a minimum is, in general, not attained; we may only construct
elements of the set arbitrarily close to the minimum. In the proof
above, selecting the point $P$ and the line $l$, for only an approximation
$e$ to the minimum $d$, would disturb the measurements, since then
$h$, while less than $e$, would not be known to be less than $d$. 

The final constructive obstacle concerns the situation in which three
points of the set $\mathscr{S}$ are assumed to lie on the selected
line $l$; the proof must decide on which side of the foot $E$ each
point lies. This amounts to deciding on which side of zero an arbitrary
real number lies, and this is precisely the nonconstructive omniscience
principle LLPO. 

In adapting Kelly's proof, these constructive obstacles will be overcome
by adding an hypothesis, by selecting a sufficiently close approximation
to the minimum distance, and by using constructive properties of the
real numbers. \\

\noindent \begin{center}
\textbf{6. Constructive version}\\
\par\end{center}

Since the Sylvester-Gallai Theorem is constructively invalid, a constructive
version must be restricted to a set of points with additional structure.
A set is \emph{discrete} if any two elements are either equal or distinct.
A set of points $\mathscr{S}$ is \emph{linearly discrete} if for
any point $P$ in $\mathscr{S}$, and for any line $l$ connecting
two distinct points of $\mathscr{S}$, either $P$ lies on $l$, or
$P$ lies outside $l$. For any set of points, these two conditions
follow from LEM; thus the constructive version of the theorem will
be classically equivalent to the traditional version. Although we
require both additional conditions, it will suffice to specify only
the second, as the lemma below will demonstrate. 

A strong definition will be used for a \emph{noncollinear} set of
points: the set contains distinct points $P,Q,R$ such that $P$ lies
outside the line $QR$. The definition of \emph{finite} set will be
as given in Section 2. 
\begin{lem*}
If a noncollinear set $\mathscr{S}$ of points on the real plane $\mathbb{R}^{2}$
is linearly discrete, then it is discrete. 
\end{lem*}
\begin{proof}
Let $A$ and $B$ be any points of $\mathscr{S}$. 

Select three distinct, noncollinear points $C_{1},C_{2},C_{3}$ of
$\mathscr{S}$. By cotransitivity, one of these, call it $P$, is
distinct from $A$. The point $B$ either lies outside the line $AP$
or it lies on $AP$. In the first case, we have $B\neq A$, so there
remains only the case in which $B\in AP$. 

Each of the points $C_{i}$ either lies on the line $AP$ or lies
outside $AP$. These noncollinear points cannot all lie on $AP$,
so there exists one, call it $Q$, with $Q\notin AP$; thus $Q\neq B$.
The point $A$ either lies outside the line $BQ$ or it lies on $BQ$.
In the first case, we have $A\neq B$ again, so only the case in which
$A\in BQ$ remains. 

Since $Q\notin AP$, the lines $AP$ and $BQ$ are distinct. The points
$A$ and $B$ are both common to each of these lines; hence $A=B$. 

The conclusion in the following constructive version of the Sylvester-Gallai
Theorem is stronger than that usually seen in classical versions;
it gives the result in an affirmative form. Rather than merely showing
that it is impossible for a third point of the given set to lie on
the selected line, the proof shows that any point of the set that
lies on the selected line must be identical with either one or the
other of the two selected points. 
\end{proof}
\begin{thm*}
Let $\mathscr{S}$ be any finite, linearly discrete, noncollinear
set of points on the real plane $\mathbb{R}^{2}$. There exist distinct
points $A,B$ in $\mathscr{S}$ such that the line $l=AB$ passes
through only these two points of $\mathscr{S}$; for any point $X$
of $\mathscr{S}$ that lies on $l$, either $X=A$ or $X=B$. 
\end{thm*}
\begin{proof}
(i) The lemma shows that the family $\mathscr{S}$ is discrete; thus
the family $\mathscr{P}$, of all pairs $(Q,R)$ of distinct points
in $\mathscr{S}$, is finite. The condition \emph{linearly discrete}
then ensures that, for each pair $(Q,R)$ in $\mathscr{P}$, the set
of points $P$ in $\mathscr{S}$, with $P$ lying outside the line
$QR$, is also finite. Thus, from the family of all ordered triads
$(P,Q,R)$ of points of $\mathscr{S}$, we may select and list those
triads with both properties, that $Q$ and $R$ are distinct, and
that $P$ lies outside the line $QR$. Hence the set of triads 
\[
\mathscr{T}=\{(P,Q,R)\,:\,P,Q,R\in\mathscr{S},\,\,Q\ne R,\,\,P\notin QR\}
\]
is finite, and the set 

\noindent 
\[
\mathscr{D}=\{\rho(P,QR)\,:\,(P,Q,R)\in\mathscr{T}\}
\]
of positive distances is also finite. 

Define $d=\min\mathscr{D}$, and let $D$ be the diameter of $\mathscr{S}$.
Since $\mathscr{D}$ is finite, the minimum exists constructively,
but the counterexample in Section 5 shows that this minimum need not
be attained by an element of $\mathscr{D}$, although close approximations
may be found. Select a real number $e$ in $\mathscr{D}$, where $e=\rho(K,AB)$
and $(K,A,B)\in\mathscr{T}$, such that $e<d\sqrt{1+d^{2}/D^{2}}$.
Define $l=AB$, and let $F$ denote the foot of the perpendicular
dropped from the point $K$ onto the line $l$. 

(ii) Coördinatize the line $l$ so that the point $F$ has coördinate
$0$, and the metric $\rho$ of $\mathbb{R}^{2}$ is preserved. Let
the coördinates of the points $A$ and $B$ be denoted $a$ and $b$.
Since $A\neq B$, we have $a\neq b$. By cotransitivity, either $a\neq0$
or $b\neq0$. By symmetry, it suffices to consider the first case.
Reversing the coördinatization if needed, we may assume that $a>0$. 

(iii) \emph{If the points $Y$ and $Z$ of $\mathscr{S}$ lie on the
line $l$ with coördinates $y$ and $z$ such that $y\geq0$ and $z\geq0$,
then $Y=Z$.}

To prove this, suppose that $y\ne z$; it suffices to consider the
case in which $z<y$. Denote by $G$ the foot of the perpendicular
dropped from the point $Z$ onto the line $YK$, and define $h=\rho(Z,G)$.
Considering the similar right triangles $YZG$ and $YKF$, we have 

\[
h/y\leq h/(y-z)=e/\sqrt{y^{2}+e^{2}}
\]

\noindent (The inequality relating the extreme terms may also be obtained
by comparing the areas of the triangles $YZK$ and $YFK$.) Since
$d\leq e$ and $y\leq\rho(Y,K)\leq D$, it follows that 
\[
h\leq e/\sqrt{1+e^{2}/y^{2}}\leq e/\sqrt{1+d^{2}/D^{2}}<d
\]

\noindent However, since $(Z,Y,K)$ belongs to the set of triads $\mathscr{T}$,
we have $h\in\mathscr{D}$, so $d\leq h$, a contradiction. Thus $y=z$. 

(iv) Suppose that $b>0$. By (iii), we then have $B=A$, a contradiction.
Thus $b\leq0$. 

(v) Now let $X$ be any point of $\mathscr{S}$ that lies on $l$,
with coördinate $x$. By cotransitivity, either $X\neq A$ or $X\neq B$.
In the first case, suppose that $x>0$. Then it follows from (iii)
that $X=A$, a contradiction. This shows that $x\leq0$, and thus,
by (iii)(reversed), we have $X=B$. Similarly, in the second case
we find that $X=A$. 
\end{proof}
There are other classical versions of the Sylvester-Gallai Theorem,
and other proofs, which might be constructivized; see  {[}KelMos58,
Wil68, Cha70, Lin88, BorMos90, Chv04, Pam09{]}. \\
\\
\textbf{Acknowledgments.} The author is grateful for useful suggestions
from the referee, and from the editor. \\

\noindent \begin{center}
\textbf{References}\\
\par\end{center}

\noindent {\small{}{[}B67{]} E. Bishop, }\emph{\small{}Foundations
of Constructive Analysis,}{\small{} McGraw-Hill, New York, 1967. MR0221878}{\small\par}

\noindent {\small{}{[}B73{]} E. Bishop, }\emph{\small{}Schizophrenia
in Contemporary Mathematics}{\small{}, AMS Colloquium Lectures, Missoula,
Montana, 1973. Reprinted in }\emph{\small{}Contemporary Mathematics}{\small{}
39:l-32, 1985. MR0788163 }{\small\par}

\noindent {\small{}{[}BB85{]} E. Bishop and D. Bridges, }\emph{\small{}Constructive
Analysis,}{\small{} Springer-Verlag, Berlin, 1985. MR0804042}{\small\par}

\noindent {\small{}{[}BR87{]} D. Bridges and F. Richman, }\emph{\small{}Varieties
of Constructive Mathematics,}{\small{} Cambridge University Press,
Cambridge, UK, 1987. MR0890955}{\small\par}

\noindent {\small{}{[}BV06{]} D. Bridges and L. Vî\c{t}\u{a},}\textsc{\small{}
}\emph{\small{}Techniques of Constructive Analysis,}{\small{} Springer,
New York, 2006. MR2253074}{\small\par}

\noindent {\small{}{[}BorMos90{]} P. Borwein and W. O. J. Moser, A
survey of Sylvester\textquoteright s problem and its generalizations,
}\emph{\small{}Aequationes Math.}{\small{} 40:111\textendash 135,
1990. MR1069788 }{\small\par}

\noindent {\small{}{[}BriMin84{]} D. Bridges and R. Mines, What is
constructive mathematics?, }\emph{\small{}Math. Intelligencer}{\small{}
6:32\textendash 38, 1984. MR0762057 }{\small\par}

\noindent {\small{}{[}Bro08{]} L. E. J. Brouwer, De onbetrouwbaarheid
der logische principes, }\emph{\small{}Tijdschrift voor Wijsbegeerte}{\small{}
2:152-158, 1908. English translation, \textquotedblleft The Unreliability
of the Logical Principles\textquotedblright , pp. 107\textendash 111
in A. Heyting (ed.), }\emph{\small{}L. E. J. Brouwer: Collected Works
1: Philosophy and Foundations of Mathematics,}{\small{} Elsevier,
Amsterdam-New York, 1975. MR0532661}{\small\par}

\noindent {\small{}{[}Cha70{]} G. D. Chakerian, Sylvester\textquoteright s
problem on collinear points and a relative, }\emph{\small{}Amer. Math.
Monthly}{\small{} 77:164\textendash 167, 1970. MR0258659}{\small\par}

\noindent {\small{}{[}Chv04{]} V. Chvátal, Sylvester-Gallai theorem
and metric betweenness,}\emph{\small{} Discrete Comput. Geom.}{\small{}
31:175\textendash 195, 2004. MR2060634}{\small\par}

\noindent {\small{}{[}Cox48{]} H. S. M. Coxeter, A problem of collinear
points,}\emph{\small{} Amer. Math. Monthly}{\small{} 55:26\textendash 28,
1948. MR0024137}{\small\par}

\noindent {\small{}{[}Cox61{]} H. S. M. Coxeter, }\emph{\small{}Introduction
to Geometry,}{\small{} Wiley, New York, 1961. MR0123930 }{\small\par}

\noindent {\small{}{[}Erd43{]} P. Erd\H{o}s, Problem for solution
4065, }\emph{\small{}Amer. Math. Monthly}{\small{} 50:65, 1943. }{\small\par}

\noindent {\small{}{[}Erd82{]} P. Erd\H{o}s, Personal reminiscences
and remarks on the mathematical work of Tibor Gallai, }\emph{\small{}Combinatorica}{\small{}
2:207\textendash 212, 1982. MR0698647}{\small\par}

\noindent {\small{}{[}KelMos58{]} L. M. Kelly and W. O. J. Moser,
On the number of ordinary lines determined by }\emph{\small{}n}{\small{}
points, }\emph{\small{}Canad. J. Math.}{\small{} 10:210\textendash 219,
1958. MR0097014}{\small\par}

\noindent {\small{}{[}Lin88{]} X. B. Lin, Another brief proof of the
Sylvester theorem, }\emph{\small{}Amer. Math. Monthly}{\small{} 95:932\textendash 933,
1988. MR979138}{\small\par}

\noindent {\small{}{[}M85{]} M. Mandelkern, Constructive mathematics,
}\emph{\small{}Math. Mag. }{\small{}58:272-280, 1985. MR0810148}{\small\par}

\noindent {\small{}{[}M89{]} M. Mandelkern, Brouwerian counterexamples,
}\emph{\small{}Math. Mag.}{\small{} 62:3-27, 1989. MR0986618}{\small\par}

\noindent {\small{}{[}Mel41{]} E. Melchior, Über Vielseite der projektiven
Ebene, }\emph{\small{}Deutsche Math.}{\small{} 5:461-475, 1941. MR0004476}{\small\par}

\noindent {\small{}{[}Pam09{]} V. Pambuccian, A reverse analysis of
the Sylvester-Gallai theorem, }\emph{\small{}Notre Dame J. Form. Log.}{\small{}
50:245\textendash 260, 2009. MR2572973 }{\small\par}

\noindent {\small{}{[}Pla05{]} J. von Plato, A constructive approach
to Sylvester\textquoteright s conjecture, }\emph{\small{}J.UCS }{\small{}11:2165\textendash 2178,
2005. MR2210695 }{\small\par}

\noindent {\small{}{[}R82{]} F. Richman, Meaning and information in
constructive mathematics, }\emph{\small{}Amer. Math. Monthly}{\small{}
89:385-388, 1982. MR0660918}{\small\par}

\noindent {\small{}{[}R99{]} F. Richman, Existence proofs, }\emph{\small{}Amer.
Math. Monthly}{\small{} 106:303-308, 1999. MR1682389}{\small\par}

\noindent {\small{}{[}R02{]} F. Richman, Omniscience principles and
functions of bounded variation, }\emph{\small{}MLQ Math. Log. Q.}{\small{}
42:111-116, 2002. MR1874208}{\small\par}

\noindent {\small{}{[}Ste44{]} R. Steinberg, Three point collinearity,
}\emph{\small{}Amer. Math. Monthly}{\small{} 51:169-171, 1944. }{\small\par}

\noindent {\small{}{[}Syl93{]} J. J. Sylvester, Mathematical question
11851, }\emph{\small{}Educational Times}{\small{} 59:98, 1893.}{\small\par}

\noindent {\small{}{[}Wil68{]} V. C. Williams, A proof of Sylvester\textquoteright s
theorem on collinear points, }\emph{\small{}Amer. Math. Monthly}{\small{}
75:980\textendash 982, 1968. MR1535105}\\

\noindent {\small{}April 7, 2015; revised March 15, 2016. }\\
Acta Mathematica Hungarica{\small{}, 150 (2016), 121\textendash 130. }{\small\par}
\end{document}